\documentclass[graybox]{svmult}

\usepackage{graphicx}           
\usepackage{newtxtext}          
\usepackage{xcolor}             
\usepackage{bm}                 
\usepackage{amsmath}            
\usepackage{amsfonts}           
\usepackage{soul}               
\usepackage{subcaption}         
\usepackage{fix-cm}             

\graphicspath{{figures/}}
\newcommand{\dx}{\Delta x}
\newcommand{\dxi}{\Delta \xi}
\newcommand\eref[1]{(\ref{#1})}
\newcommand\fref[1]{Fig.~\ref{#1}}

\begin{document}

\title*{Spline-Based Stochastic Collocation Methods for Uncertainty Quantification in Nonlinear Hyperbolic PDEs}
\titlerunning{Spline-Based Stochastic Collocation Methods for Uncertainty Quantification}
\author{Alina Chertock\orcidID{0000-0003-4978-1314} and\\
Arsen S. Iskhakov\orcidID{0000-0003-4064-9894} and\\
Safa Janajra\orcidID{0000-0002-6564-1971} and\\
Alexander Kurganov\orcidID{0000-0002-2263-8772}}
\authorrunning{A. Chertock, A. S. Iskhakov, S. Janajra and A. Kurganov} 
\institute{Alina Chertock, Safa Janajra \at Department of Mathematics, North Carolina State University, Raleigh, NC, 
USA \\ \email{chertock@math.ncsu.edu, sjanajr@ncsu.edu}
\and Arsen S. Iskhakov \at Alan Levin Department of Mechanical and Nuclear Engineering, Kansas State University, Manhattan, KS, USA\\
\email{aiskhak@ksu.edu}
\and Alexander Kurganov \at Department of Mathematics, Shenzhen International Center for Mathematics, and Guangdong Provincial Key 
Laboratory of Computational Science and Material Design, Southern University of Science and Technology, Shenzhen, China\\ 
\email{alexander@sustech.edu.cn}}

\maketitle

\abstract{In this paper, we study the stochastic collocation (SC) methods for uncertainty quantification (UQ) in hyperbolic systems of 
nonlinear partial differential equations (PDEs). In these methods, the underlying PDEs are numerically solved at a set of collocation points
in random space. A standard SC approach is based on a generalized polynomial chaos (gPC) expansion, which relies on choosing
the collocation points based on the prescribed probability distribution and approximating the computed solution by a linear combination of
orthogonal polynomials in the random variable. We demonstrate that this approach struggles to accurately capture discontinuous solutions, 
often leading to oscillations (Gibbs phenomenon) that deviate significantly from the physical solutions. We explore alternative SC methods,
in which one can choose an arbitrary set of collocation points and employ shape-preserving splines to interpolate the solution in a random
space. Our study demonstrates the effectiveness of spline-based collocation in accurately capturing and assessing uncertainties while
suppressing oscillations. We illustrate the superiority of the spline-based collocation on two numerical examples, including the inviscid
Burgers and shallow water equations.}

\section{Introduction}\label{sec1}
Numerous scientific problems encompass inherent uncertainties arising from a variety of factors. Within the context of partial differential 
equations (PDEs), uncertainties can be characterized using random variables. This study focuses on hyperbolic systems of conservation and 
balance laws, examining their behavior under uncertain conditions. In the one-dimensional (1-D) case, the formulation of these systems is 
expressed as
\begin{equation}
\bm U_t+\bm F(\bm U)_x=\bm S(\bm U),
\label{1}
\end{equation}
where $x$ is the spatial variable, $t$ is time, $\bm U(x,t;\xi)\in\mathbb R^m$ is an unknown vector-function, 
$\bm F(\bm U):\mathbb R^m\to\mathbb R^m$ are the flux functions, and $\bm S(\bm U)$ is a source term. Furthermore, we assume 
that $\xi\in\Xi\subset\mathbb R$ are real-valued random variable with $(\Xi,{\cal F},\nu)$ being the underlying probability space. 
Here, $\Xi$ is a set of events, ${\cal F}(\Xi)$ is the $\sigma$-algebra of Borel measurable sets, and $\nu(\xi):\Xi\to\mathbb R_+$ is the 
probability density function (PDF), $\nu\in L^1(\Xi)$.

The system \eref{1} emerges in various applications, such as fluid dynamics, geophysics, electromagnetism, meteorology, and astrophysics. 
It is crucial to assess uncertainties inherent in input quantities, as well as in the initial and boundary conditions resulting from 
empirical approximations or measurement errors. This quantification is vital for performing the sensitivity analysis and offers valuable 
insight to improve the accuracy of the studied model.

This paper focuses on developing accurate and robust numerical techniques for quantifying uncertainties in \eref{1}. Among the 
various existing methods, Monte Carlo-type methods (see, e.g., \cite{AM17,MSS13}) are reliable but computationally intensive due to the 
substantial number of realizations required to approximate statistical moments accurately. Another commonly used approach is based on 
generalized polynomial chaos (gPC) methods, in which the solution is expressed as a series of orthogonal polynomials with respect to the 
probability density in $\xi$ \cite{LK10,PIN-book}. There are two types of gPC methods: intrusive and non-intrusive. Intrusive approaches, 
such as gPC stochastic Galerkin (gPC-SG) methods, substitute gPC expansions into the governing equations. These expansions are then 
projected using a Galerkin approximation to derive a system of deterministic PDEs for the expansion coefficients; see, e.g., 
\cite{TLNE10a,TLNE10}. Solving these coefficient equations provides the statistical moments of the original solution of the uncertain 
problem. On the other hand, non-intrusive algorithms, such as gPC stochastic collocation (gPC-SC) methods, aim to satisfy the governing 
equations at discrete nodes, called collocation points, in the random space. They employ a deterministic numerical solver, utilizing 
interpolation and quadrature rules to numerically evaluate the PDF and/or statistical moments \cite{Xiu09,SMS13}.

The application of gPC-SG and gPC-SC methods to nonlinear hyperbolic systems \eref{1} poses several challenges. Although spectral-type 
gPC-based methods exhibit rapid convergence for solutions that depend smoothly on random parameters, a significant issue arises when 
solutions contain shock waves and other nonsmooth structures, which is a generic case for nonlinear hyperbolic PDEs (even if initial data 
are smooth). Despite the discontinuities manifesting in the spatial variable, their propagation speed can be influenced by uncertainty, 
introducing discontinuities in the random variable and causing Gibbs-type phenomena \cite{LKNG04}. Another unresolved matter involves 
representing strictly positive quantities, such as water depth in shallow water equations, and imposing discrete bound-preserving 
constraints \cite{DEN21a,GHS19,SS18}. 

In this paper, we concentrate on alternative spline-based stochastic collocation (SC) methods, which are positivity-preserving and do not
suffer from Gibbs-type oscillations. The proposed methods utilize a solution obtained by a deterministic numerical solver implemented
repeatedly on an arbitrarily selected set of collocation points in the random variable $\xi$. Specifically, at each collocation point, we
solve the corresponding deterministic PDEs by a semi-discrete second-order central-upwind scheme from \cite{KNP,KPshw}. As a result, at the
final computational time, for each discrete value of the spatial variable $x$, we obtain an approximation of $\bm U$ as a discrete function
of $\xi$. Equipped with these point values, we employ spline-based interpolations in random space and use the obtained global (in $\xi$)
solution to calculate the stochastic moments. In order to enforce the non-oscillatory and positivity-preserving properties of the
interpolated solution, we use shape-preserving (SP) rational quartic splines from \cite{ZH15}. We conduct a comparative numerical study of
the gPC-based and spline-based SC approaches to quantify uncertainties in \eref{1}. Our findings demonstrate the superior efficacy of the
proposed spline-based SC methods when applied to the inviscid Burgers and shallow water equations.

\section{Methodology}\label{sec2}
In this section, we describe the gPC- and spline-based SC approaches applied to \eref{1}. 

We start by selecting a set of collocation points $\xi_\ell,\ \ell=1,\ldots,L$ and numerically solving the following deterministic systems:
\begin{equation}
\bm U_t(x,t;\xi_\ell)+\bm F(\bm U(x,t;\xi_\ell))_x=\bm S(\bm U(x,t;\xi_\ell)),\quad\ell=1,\ldots,L,
\label{2}
\end{equation}
until the final time $T$ (one can use one's favorite numerical method for solving \eref{2}). Then, for each discrete node of the spatial 
variable denoted below by $\tilde x$, we use either the gPC (\S\ref{sec2.1}) or spline (\S\ref{sec2.2}) interpolation to approximate the
numerical solution and its stochastic moments, that is, the mean and variance or standard deviation for each component $U$ of $\bm U$:
\begin{equation}
\begin{aligned}
&\mathbb E_\xi[U]:=\int\limits_\Xi U(\tilde x,T;\xi)\nu(\xi){\rm d}\xi,\\
&{\rm Var}[U]:=\mathbb E_\xi[U^2]-\mathbb E_\xi[U]^2,\quad\sigma[U]:=\sqrt{{\rm Var}[U]}.
\end{aligned}
\label{3}
\end{equation}

\subsection{gPC Interpolation}\label{sec2.1}
The gPC interpolation in random space represents the solution as a generalized discrete Fourier series in terms of orthonormal polynomials, 
$\{\Phi_i(\xi)\}_{i=0}^N$, selected based on the PDF:
\begin{equation}
\bm U(\tilde x,T;\xi)\approx\bm U^N(\tilde x,T;\xi):=\sum_{i=0}^N\widehat{\bm U}_i(\tilde x,T)\,\Phi_i(\xi),
\label{4}
\end{equation}
where $\widehat{\bm U}_i(\tilde x,T)$ are deterministic Fourier coefficients. 

It is well-known that for large values of $N$, the polynomial interpolation \eref{4} may be very oscillatory. To minimize the oscillations,
one needs to choose the roots of $\Phi_{N+1}(\xi)$ as collocation points $\xi_\ell,\ \ell=1,\ldots,L$ with $L=N+1$. In this case, the
Fourier coefficients can be computed with the help of the discrete Fourier transform:
\begin{equation}
\widehat{\bm U}_i(\tilde x,T)=\sum_{\ell=1}^{N+1}\bm U_\ell\,\Phi_i(\xi_\ell)\,\omega_\ell,\quad i=0,\ldots,N,
\label{5}
\end{equation}
where $\bm U_\ell:\approx\bm U(\tilde x,T;\xi_\ell)$ and $\omega_\ell$ are the Gauss quadrature weights corresponding to the PDF $\nu(\xi)$.
These coefficients are also used to calculate the stochastic moments of for each component $U$ of the computed solution $\bm U$:
\begin{equation*}
\mathbb E_\xi[U^N]=\widehat U_0,\quad{\rm Var}[U^N]=\sum_{i=1}^N\widehat U^2_i.
\end{equation*}

It should be observed that the gPC interpolation \eref{4} is exponentially accurate for smooth solutions but suffers from the Gibbs-type
phenomenon when discontinuities appear in the numerical solution, which is a generic case for nonlinear hyperbolic PDEs. Therefore, in the 
next two sections, we turn our attention to alternative interpolation methods that lead to non-oscillatory approximation of 
$\bm U(\tilde x,T;\xi)$.

\subsection{Shape-Preserving (SP) Spline Interpolation}\label{sec2.2}
In order to suppress Gibbs-type oscillations, one can use cubic B-splines (see, e.g., \cite{Cox72,deB72,deBP03,ZM97}), which also retain the
positivity of the interpolated data, but as we demonstrate in our numerical experiments below, B-spline approximations may oversmear shock 
discontinuities. Moreover, B-splines do not necessarily maintain the monotonicity and/or convexity of the interpolated data. 

Therefore, in this paper, we use (provably) SP rational quartic interpolation splines from \cite{ZH15}, specifically designed to preserve
the shape of positive, monotonic, and convex solutions. These SP splines also ensure $C^2$ continuity of the constructed interpolant.

Let us describe the SP spline for a certain component $U$ of $\bm U$. We denote by
$\dxi_\ell:=\xi_{\ell+1}-\xi_\ell,\ \ell=1,\ldots,L-1$, and introduce the following quantities:
\begin{equation*}
\begin{aligned}
&\delta_\ell:=\frac{1}{\dxi_\ell+\dxi_{\ell-1}}\left(\frac{U_{\ell+1}-U_\ell}{\dxi_\ell}
-\frac{U_\ell-U_{\ell-1}}{\dxi_{\ell-1}}\right),\\
&A_1:=\dxi_1\delta_2,\quad A_\ell:=\dxi_\ell\delta_\ell,\\
&B_{\ell-1}:=\dxi_{\ell-1}\delta_\ell,\quad B_{L-1}:=\dxi_{L-1}\delta_{L-1},
\end{aligned}
\qquad\ell=2,\ldots,L-1.
\end{equation*}
For each subinterval $\left[\xi_\ell,\xi_{\ell+1}\right]$, $\ell=1,2,\dots,L-1$, the computed solution is interpolated by the spline given
by
\begin{equation*}
\begin{aligned}
U(\tilde x,T;\xi)&\approx S\left(\tilde x,T;\xi\right):=(1-\tau_\ell)U_\ell+\tau_\ell U_{\ell+1}\\
&-\frac{\dxi_\ell(1-\tau_\ell)\tau_\ell\left[(1-\tau_\ell)^2A_\ell+\lambda_\ell(1-\tau_\ell)\tau_\ell A_\ell B_\ell+
\tau_\ell^2B_\ell\right]}{Q_\ell},
\end{aligned}
\end{equation*}
where
\begin{equation*}
Q_\ell=1+(1-\tau_\ell)\tau_\ell\big[(1-\tau_\ell)(B_\ell\lambda_\ell+\mu_\ell)+\tau_\ell(A_\ell\lambda_\ell+\mu_{\ell+1})\big],
\end{equation*}
where $\tau_\ell:=(\xi-\xi_\ell)/{\dxi_\ell}$, and $\lambda_\ell$ and $\mu_\ell$ are two local tension shape parameters. Note that 
specific choices of $\lambda_\ell$ and $\mu_\ell$ guarantee the preservation of monotonicity, positivity, and convexity of the 
interpolation spline; see \cite{ZH15}.

Once the spline function is constructed, one can use one's favorite quadrature rule to numerically approximate the stochastic moments in 
\eref{3}.

\section{Numerical Examples}\label{sec3}
In this section, we illustrate the performance of the gPC and spline-based SC approaches on numerical examples for the inviscid Burgers and
shallow water equations. We consider a random variable $\xi$ uniformly distributed on the interval $[-1,1]$ ($\xi\in\mathcal{U}[-1,1]$),
which induces the usage of the Legendre polynomials $\Phi_i(\xi)$ in the generalized Fourier expansion \eref{4} and Gauss-Legendre
quadratures for computing the gPC coefficients in \eref{5}. For the spline interpolation, we take equally-spaced collocation nodes
$\xi_\ell$, but note that in principle, they can be chosen arbitrarily {\color{black}since, unlike the gPC interpolation, the spline one can
be efficiently conducted on any set of nodes}. Recall that the SP splines require a specific choice of the local tension parameters
$\lambda_\ell$ and $\mu_\ell$ to guarantee the shape-preserving properties of the constructed interpolant; see
\cite[Formulae (12), (19), (24), and (28)]{ZH15} for details. 

In this work, we use semi-discrete second-order central-upwind schemes from \cite{KNP,KPshw} to numerically solve the
deterministic systems \eref{2}. The semi-discretization results in the system of ODEs, which is integrated in time using the three-stage
third-order strong stability preserving (SSP) Runge-Kutta method (see, e.g., \cite{GKS,GST}) with the time step chosen according to the CFL
number 0.45. The central-upwind scheme employs the generalized minmod limiter with parameter $\theta=1.3$; see \cite{KNP,KPshw} for details.

\subsection*{Example 1---Burgers Equation}
We start with the inviscid Burgers equation,
\begin{equation}
U_t+\left(\frac{U^2}{2}\right)_x=0,
\label{10}
\end{equation}
considered subject to the following stochastic initial data:
\begin{equation}
U(x,0;\xi)=\begin{cases}
2,\quad x>0.1\xi,\\
1,\quad x<0.1\xi.
\end{cases}
\label{11}
\end{equation}
The problem is numerically solved in the physical domain $x\in[-1,1]$ with free boundary conditions. 

\begin{figure}[ht!]
\centering
\begin{subfigure}[t]{0.32\textwidth}
\centering
\includegraphics[width=1.0\textwidth]{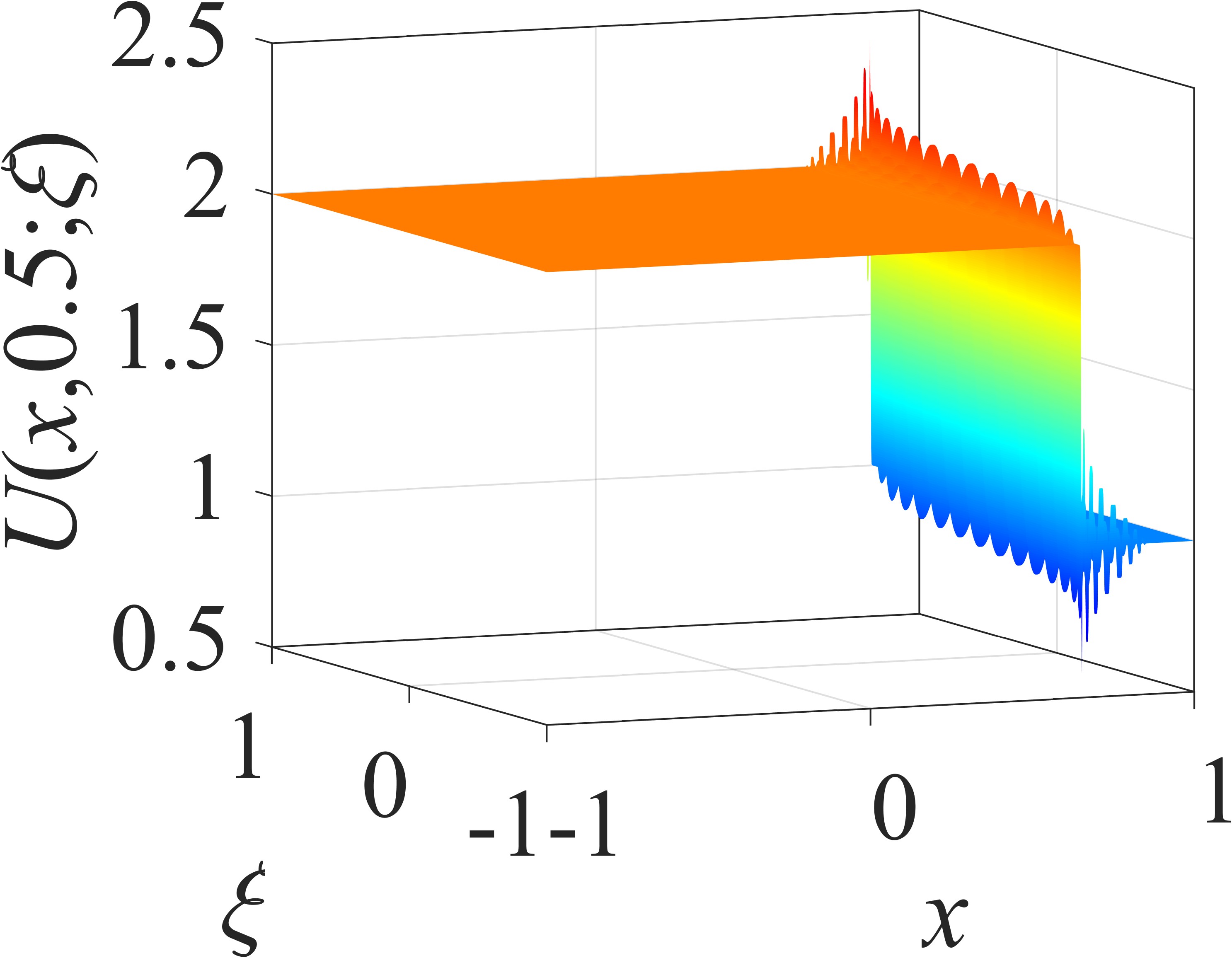}
\caption{}
\label{fig3.1a}
\end{subfigure}
\hspace{0.1cm}
\centering
\begin{subfigure}[t]{0.29\textwidth}
\centering
\includegraphics[width=1.0\textwidth]{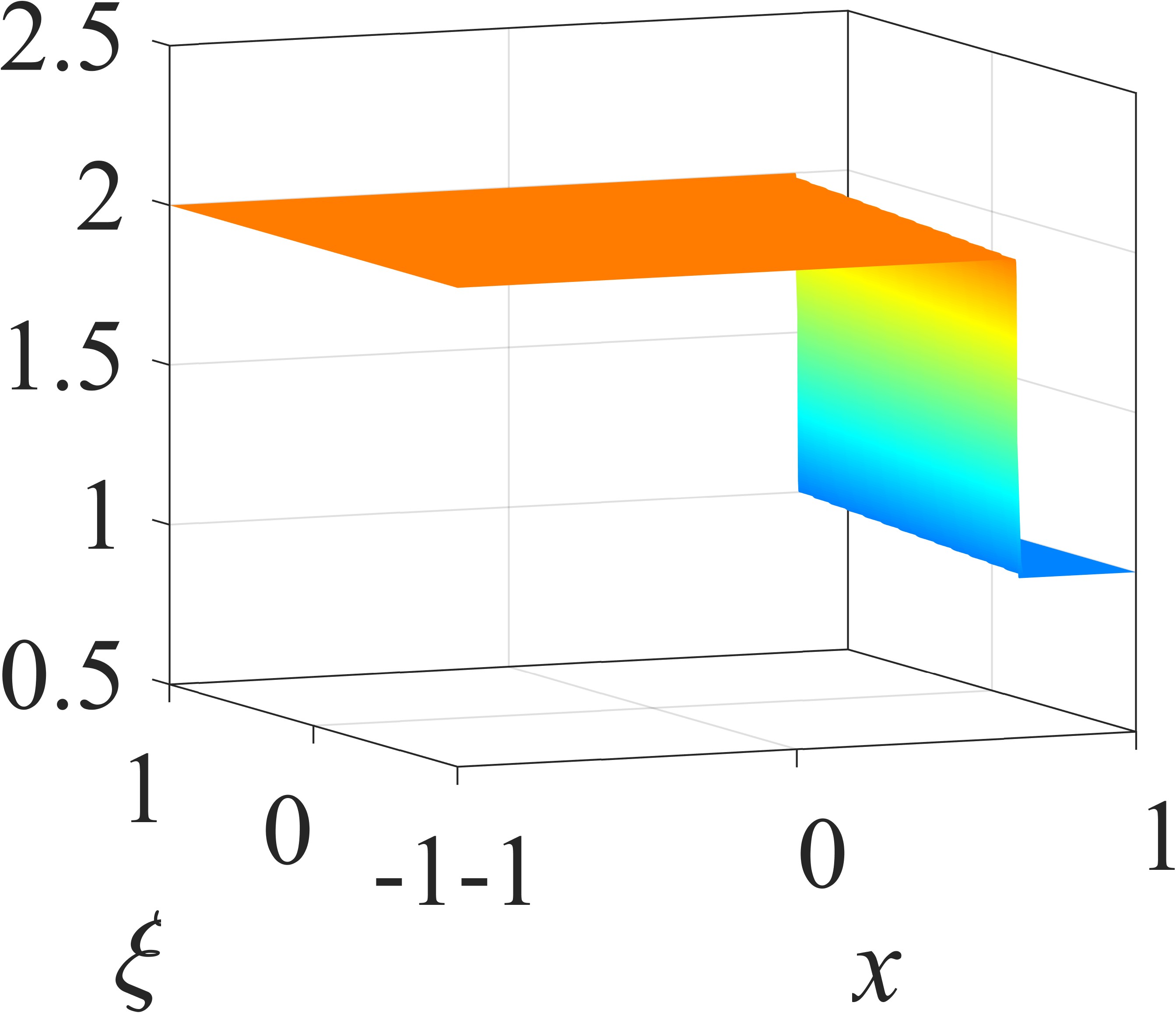}
\caption{}
\label{fig3.1b}
\end{subfigure}
\hspace{0.1cm}
\centering
\begin{subfigure}[t]{0.35\textwidth}
\centering
\includegraphics[width=1.0\textwidth]{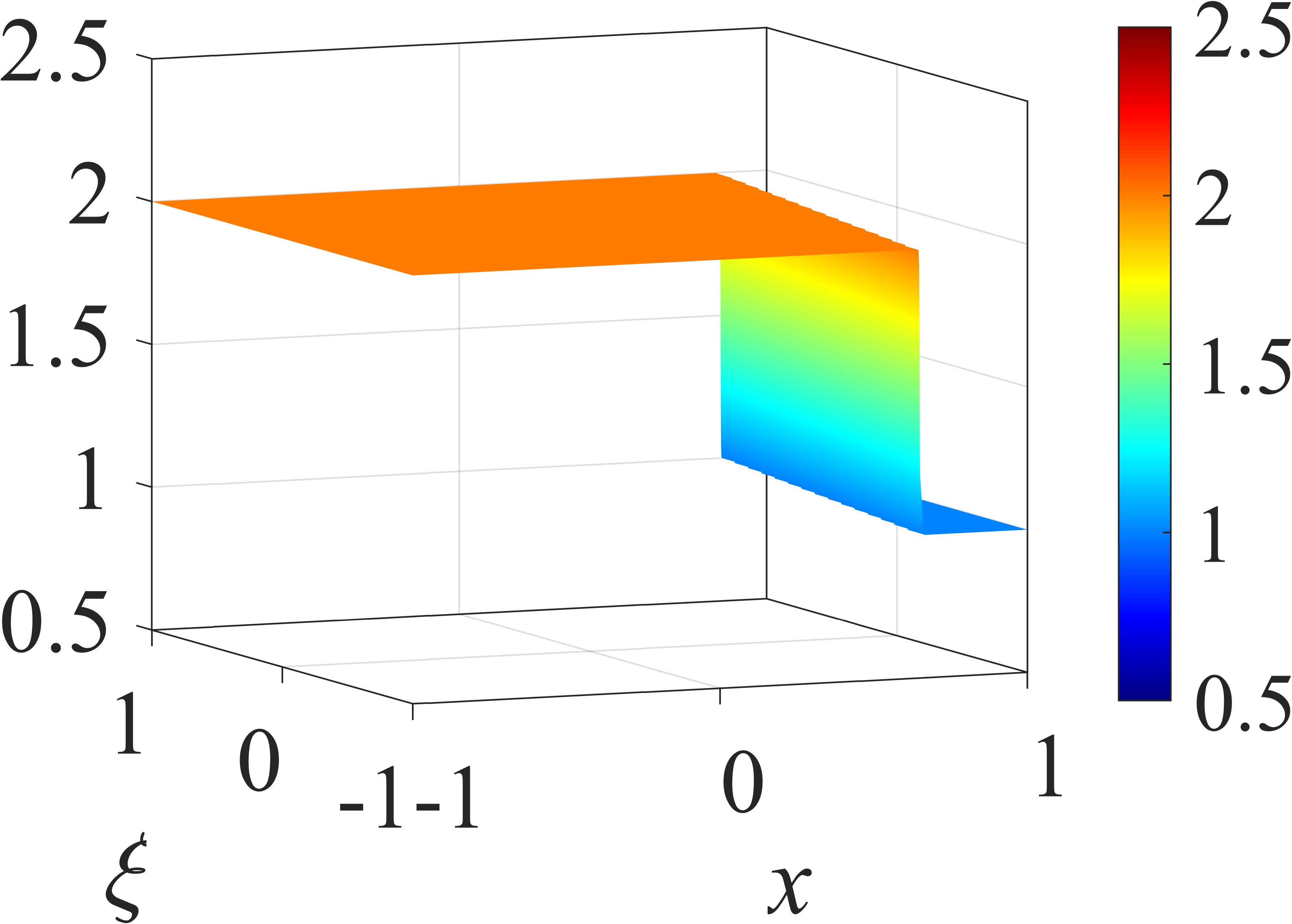}
\caption{}
\label{fig3.1c}
\end{subfigure}
\captionsetup{labelfont=bf}
\caption{Example 1: $U(x,0.5;\xi)$ obtained using (\protect\subref{fig3.1a}) gPC expansion, (\protect\subref{fig3.1b}) B-splines, and 
(\protect\subref{fig3.1c}) SP splines.}
\label{fig3.1}
\end{figure}

\begin{figure}[ht!]
\centering
\begin{subfigure}[t]{0.32\textwidth}
\centering
\includegraphics[width=1.0\textwidth]{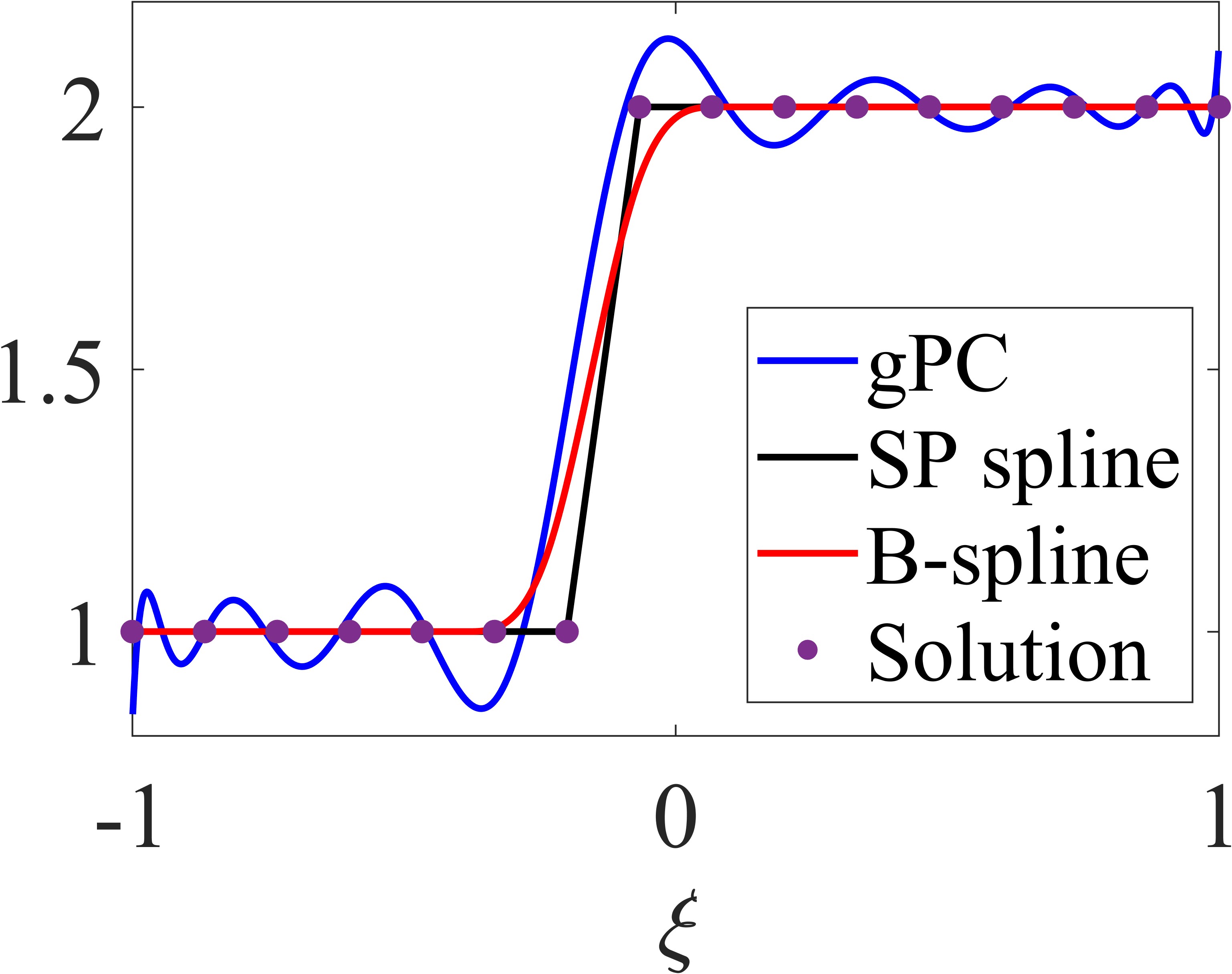}
\caption{}
\label{fig3.2a}
\end{subfigure}
\centering
\begin{subfigure}[t]{0.32\textwidth}
\centering
\includegraphics[width=1.0\textwidth]{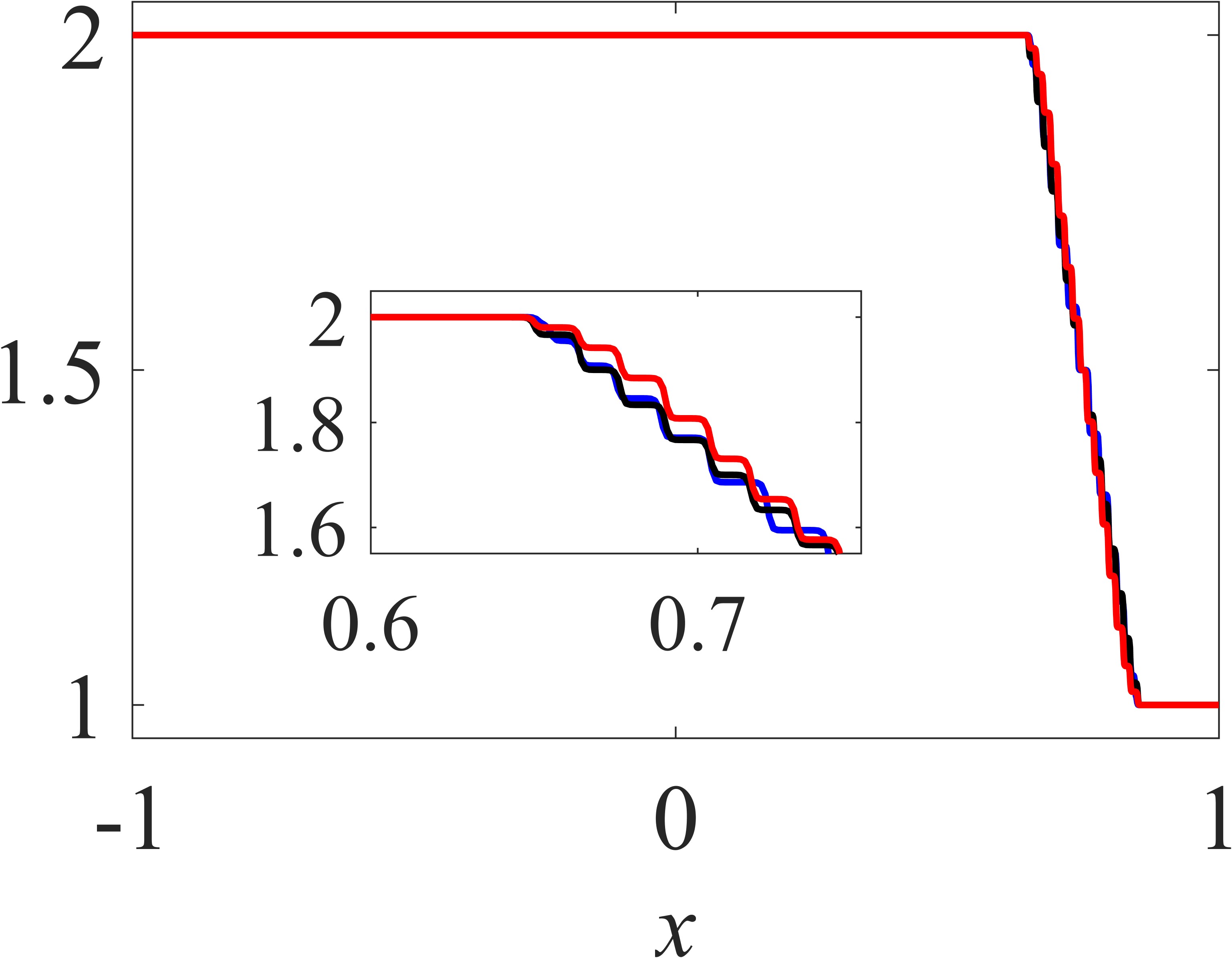}
\caption{}
\label{fig3.2b}
\end{subfigure}
\begin{subfigure}[t]{0.32\textwidth}
\centering
\includegraphics[width=1.0\textwidth]{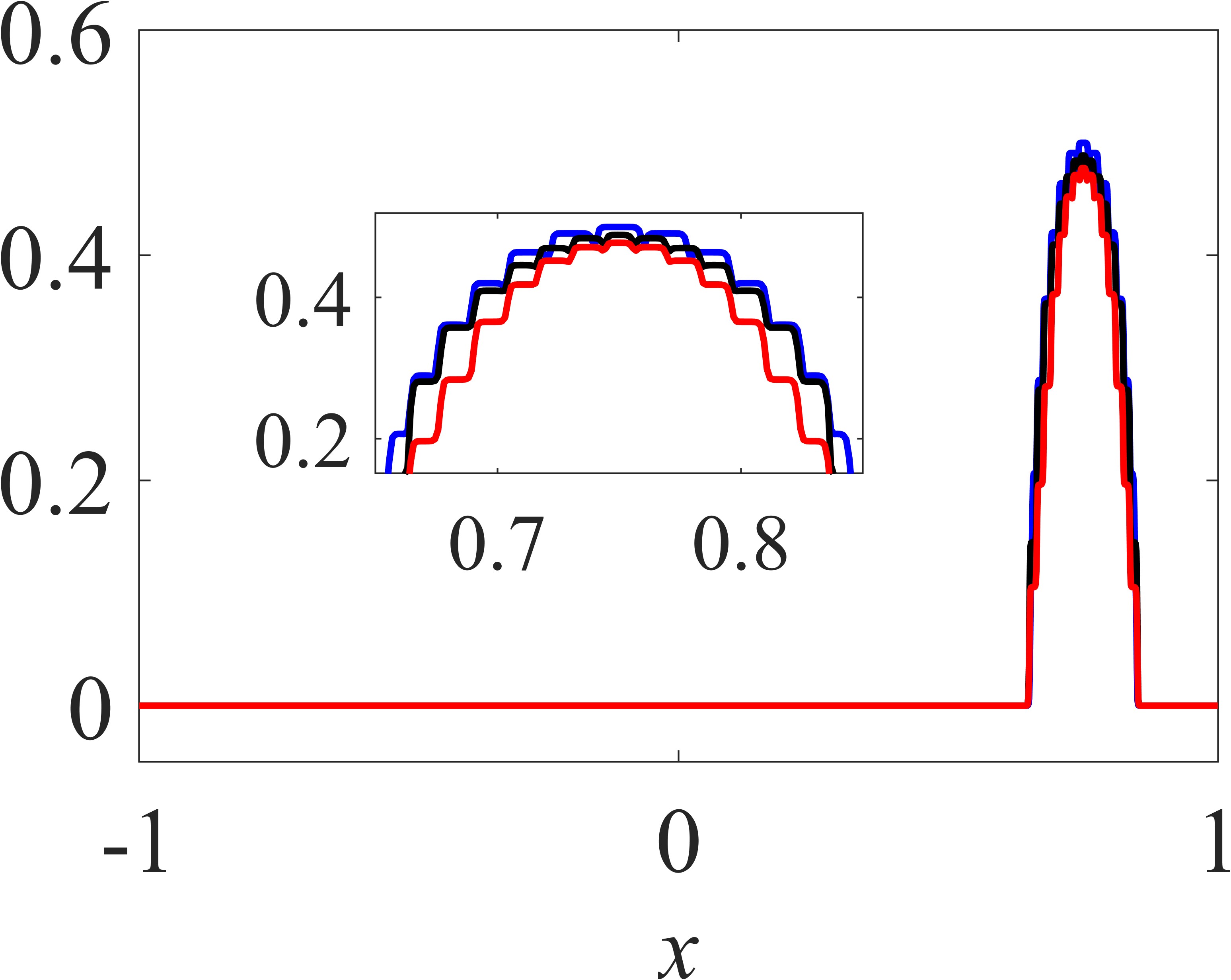}
\caption{}
\label{fig3.2c}
\end{subfigure}
\captionsetup{labelfont=bf}
\caption{Example 1: (\protect\subref{fig3.2a}) $U(0.734,0.5;\xi)$, (\protect\subref{fig3.2b}) mean, and (\protect\subref{fig3.2c}) standard 
deviation obtained using the gPC expansion and splines. {\color{black}Panels (b) and (c) contain a zoom at the areas of interest.}}
\label{fig3.2}
\end{figure}

To numerically solve the problem \eref{10}--\eref{11}, we choose $L=16$ collocation points and perform deterministic simulations on a
uniform spatial mesh consisting of cells of size $\dx=1/800$ until the final time $T=0.5$. In \fref{fig3.1}, we plot the solution
$U(x,0.5;\xi)$ obtained by interpolating the computed data in the random variable $\xi$ by the gPC and spline-based approaches. For the
latter, we present results produced by both the B-splines and SP splines. As one can see, the gPC approach exhibits Gibbs-type oscillations
near the discontinuity, whereas both spline-based approximations are oscillation-free. \fref{fig3.2a} additionally shows the profile of
$U(0.734,0.5;\xi)$ to provide more evidence that the shock propagates from the physical space to the random space. From this figure, one can
also observe that B-splines smear the discontinuity compared to the SP splines. The mean and variance values obtained with all three
approaches appear similar, as demonstrated in \fref{fig3.2b} and \ref{fig3.2c}.

\subsection*{Example 2---Shallow Water Equations}
In this example, we consider the Saint-Venant system of shallow water equations given by \eref{1} with
\begin{equation*}
\bm U=(h,hu)^\top,\quad\bm F(\bm U)=\Big(hu,hu^2+\frac{g}{2}h^2\Big)^\top,\quad\bm S=(0,-ghZ_x)^\top,
\end{equation*}
where $h(x,t;\xi)$ is the water depth, $u(x,t;\xi)$ is the velocity, $Z(x;\xi)$ is the bottom topography, and $g$ is the constant
acceleration due to gravity (we take $g=1$). 

The Saint-Venant system is considered in the physical domain $x\in[-1,1]$ subject to free boundary condition, deterministic initial data
for the water surface $w=h+Z$ and velocity $u$,
\begin{equation*}
w(x,0;\xi)=\left\{\begin{aligned}&1,&&x<0,\\&0.5,&&x>0,\end{aligned}\right.\qquad u(x,0;\xi)\equiv0,
\end{equation*}
and stochastic bottom topography
\begin{equation*}
Z(x)=\begin{cases}
0.125\xi+0.125(\cos(5\pi x)+2), &|x|<0.2,\\
0.125\xi+0.125,&\mbox{otherwise}.
\end{cases}
\end{equation*}

We numerically solve the deterministic systems \eref{2} on a uniform mesh consisting of cells of size $\dx=1/400$ until the final time
$T=0.8$ on a set of collocation points. Here, we examine the performance of the proposed methods with the number of collation points set to
either $L=16$ or $L=32$.

\begin{figure}[ht!]
\centering
\begin{subfigure}[t]{0.32\textwidth}
\centering
\includegraphics[width=1.0\textwidth]{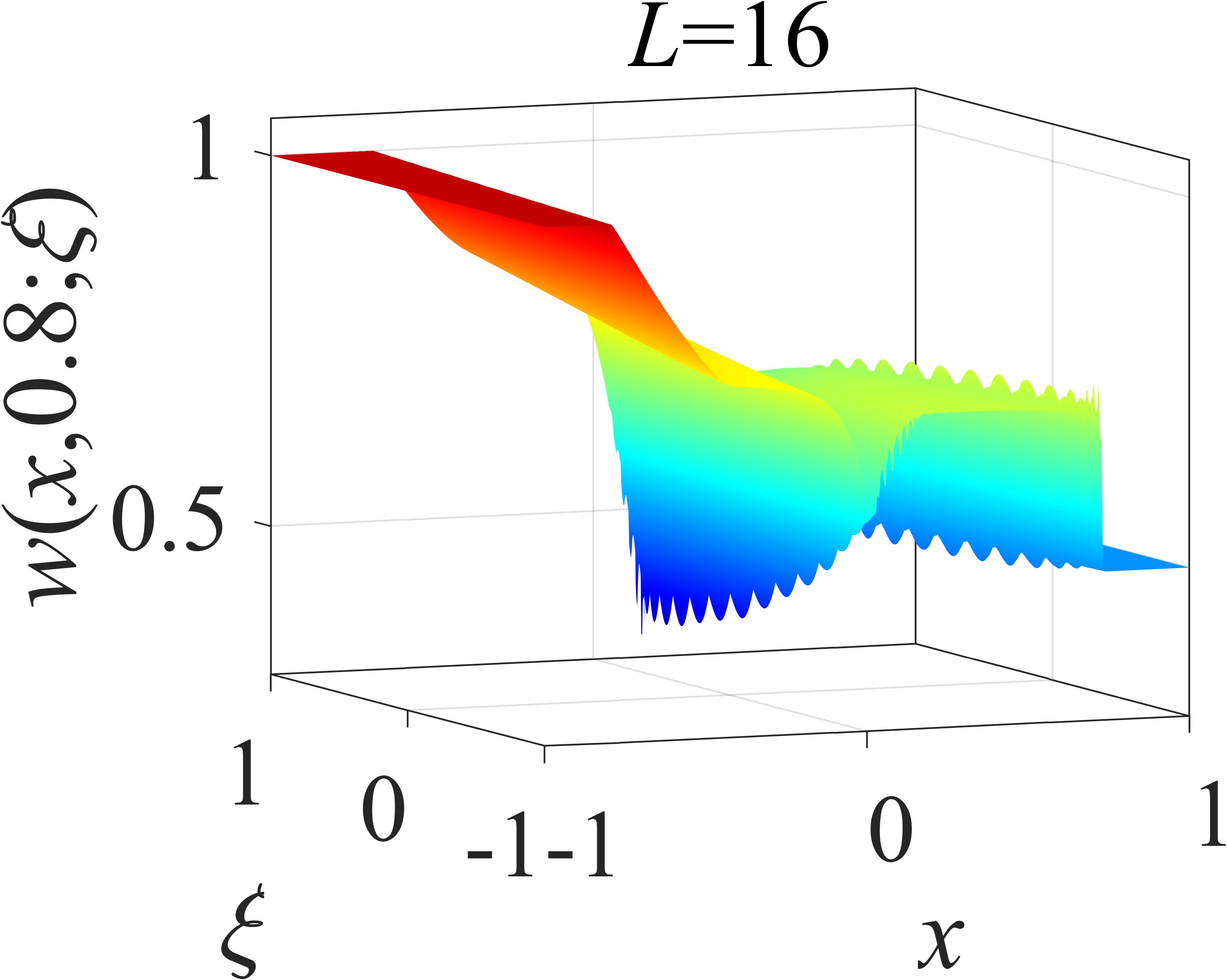}
\label{fig3.3aa}
\vspace{-0.25cm}
\end{subfigure}
\hspace{0.05cm}
\centering
\begin{subfigure}[t]{0.29\textwidth}
\centering
\includegraphics[width=1.0\textwidth]{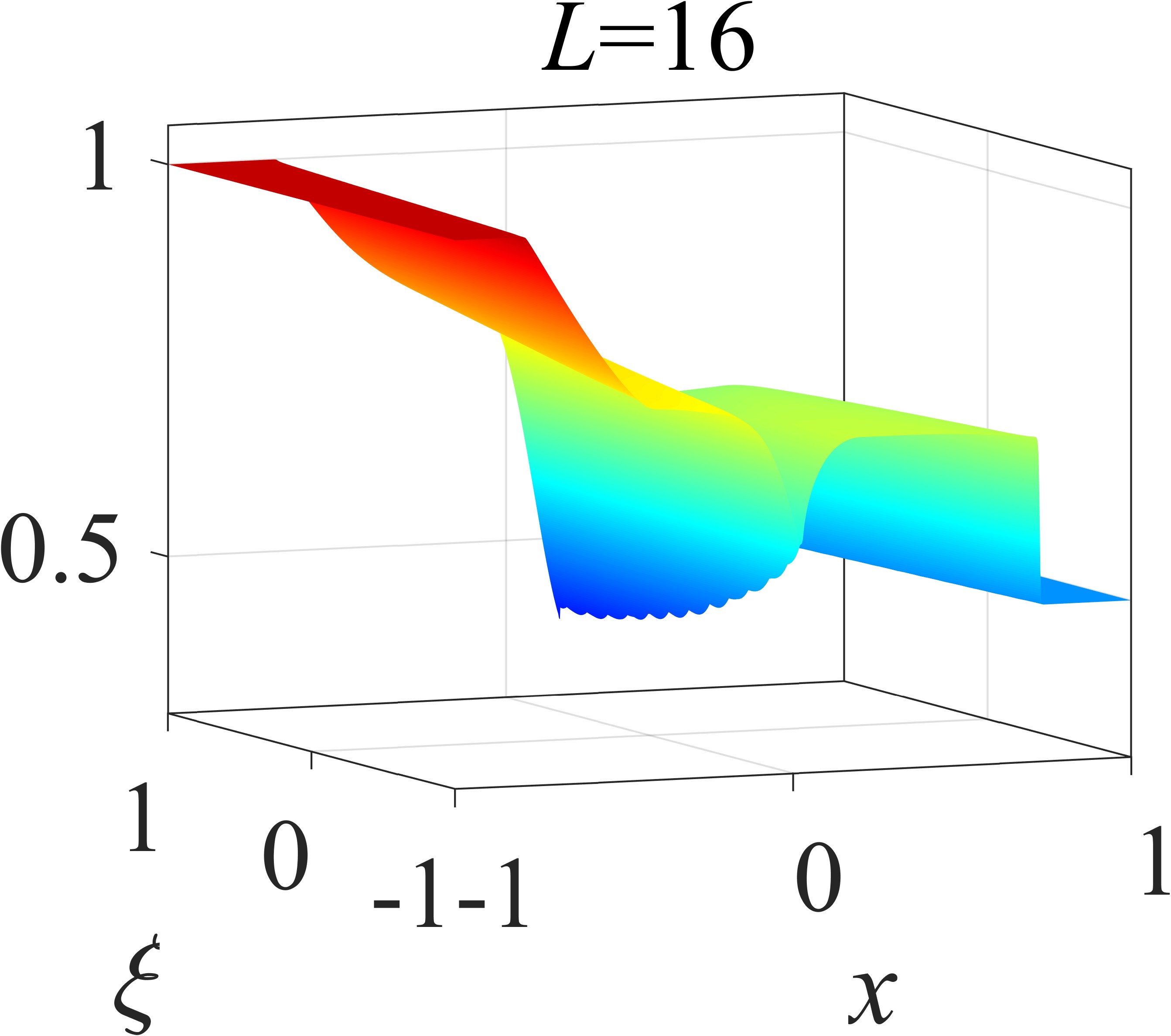}
\label{fig3.3bb}
\vspace{-0.25cm}
\end{subfigure}
\hspace{0.05cm}
\centering
\begin{subfigure}[t]{0.35\textwidth}
\centering
\includegraphics[width=1.0\textwidth]{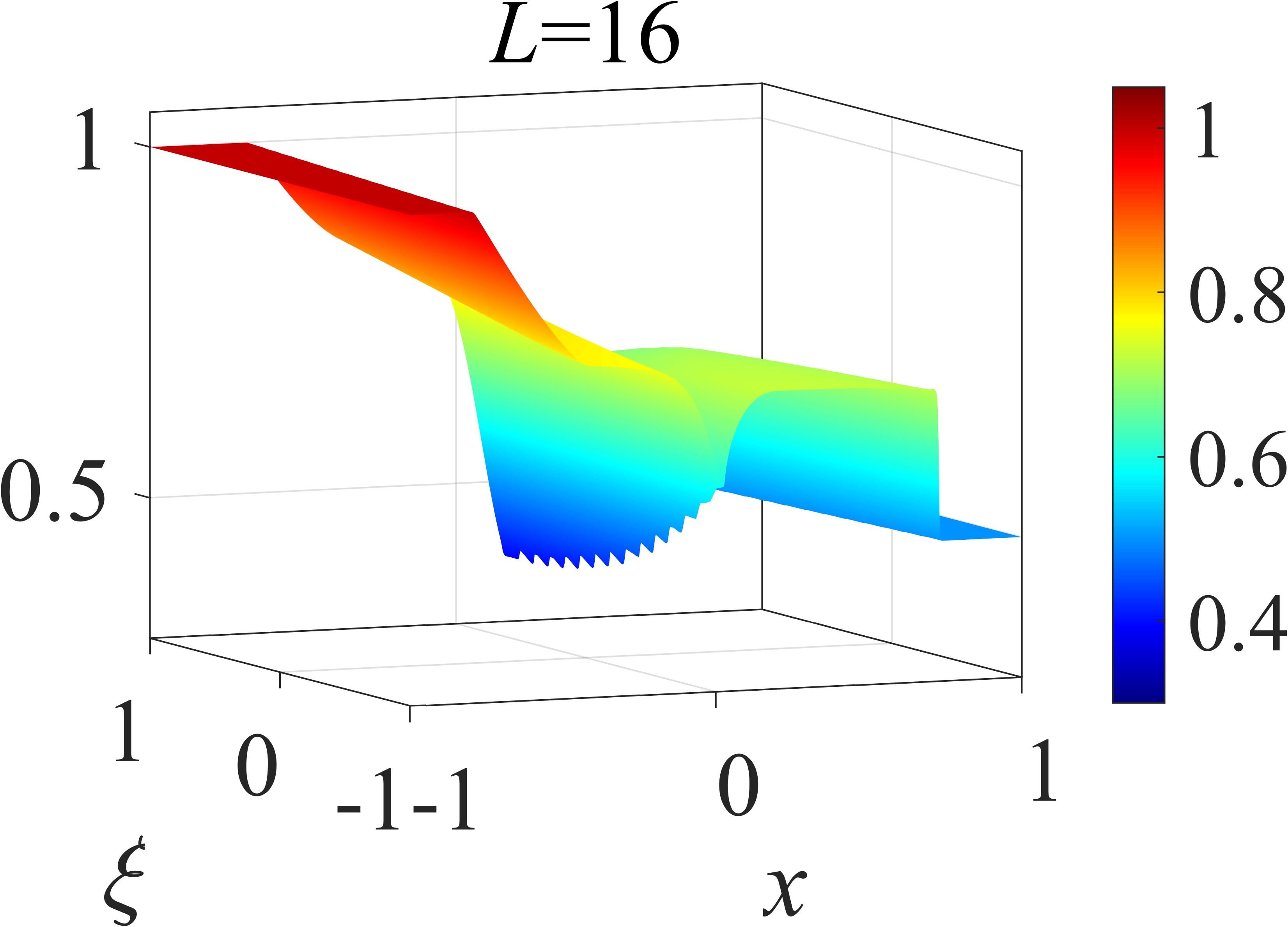}
\label{fig3.3cc}
\vspace{-0.25cm}
\end{subfigure}
\centering
\begin{subfigure}[t]{0.32\textwidth}
\centering
\includegraphics[width=1.0\textwidth]{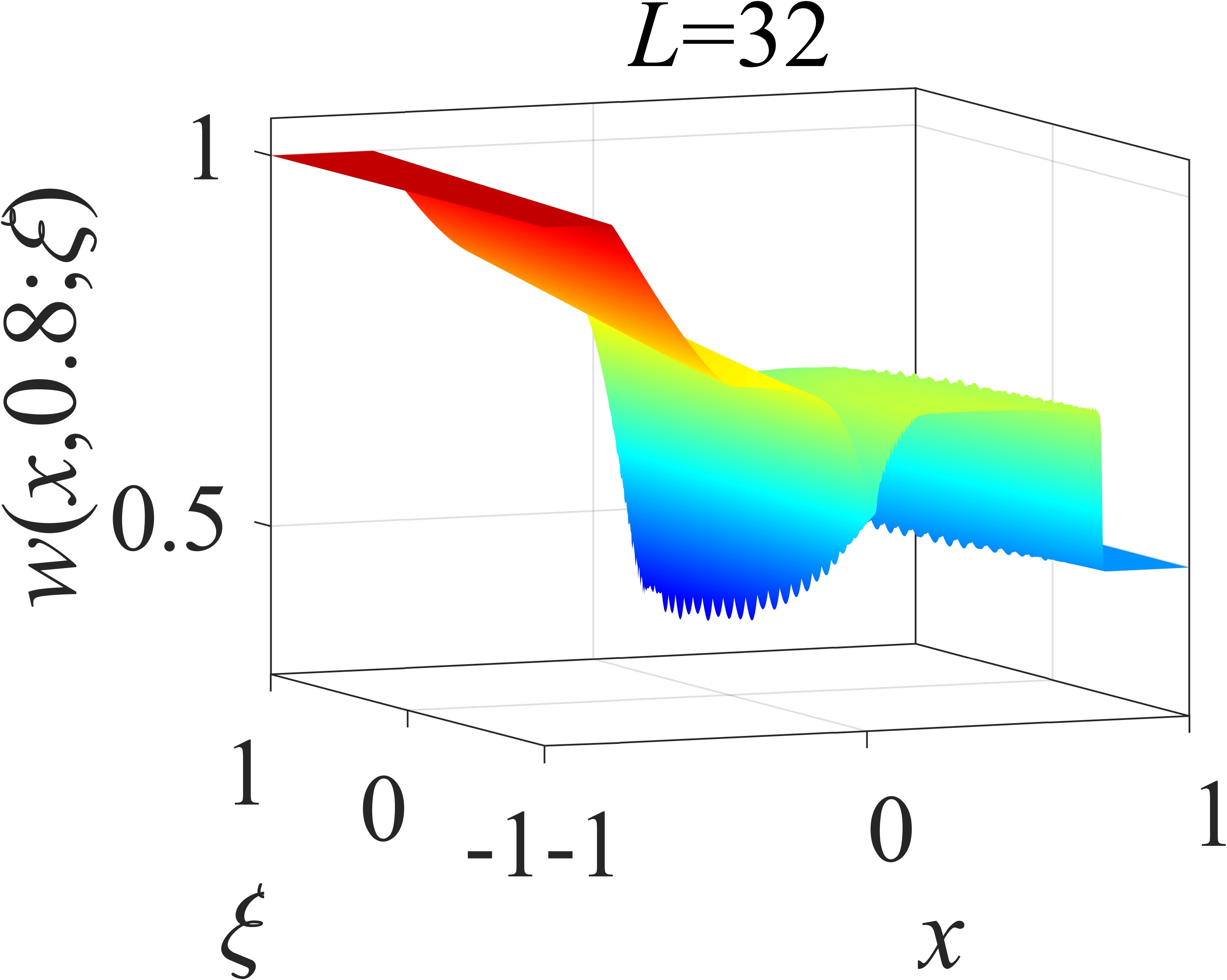}
\caption{}
\label{fig3.3a}
\end{subfigure}
\hspace{0.05cm}
\centering
\begin{subfigure}[t]{0.29\textwidth}
\centering
\includegraphics[width=1.0\textwidth]{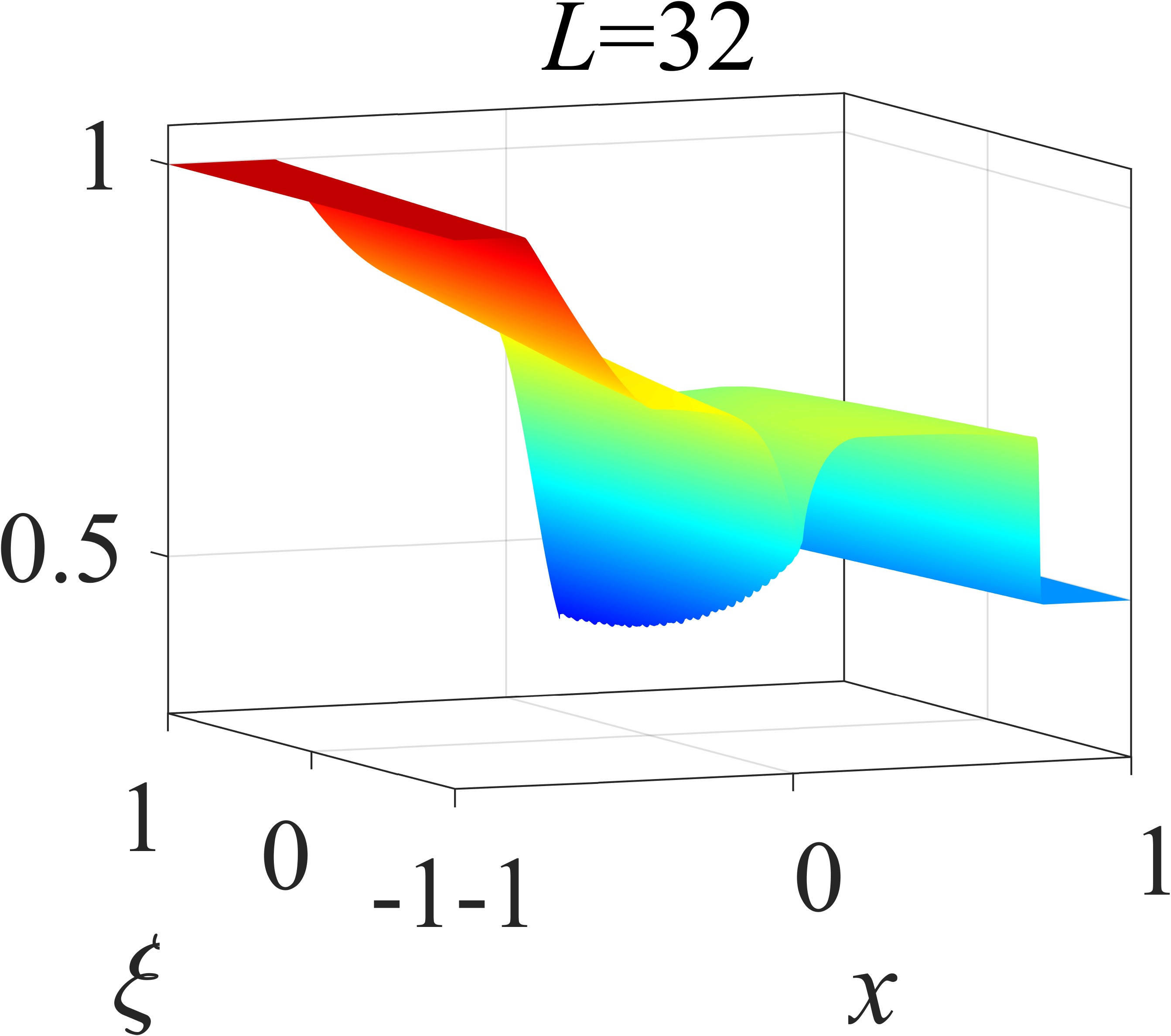}
\caption{}
\label{fig3.3b}
\end{subfigure}
\hspace{0.05cm}
\centering
\begin{subfigure}[t]{0.35\textwidth}
\centering
\includegraphics[width=1.0\textwidth]{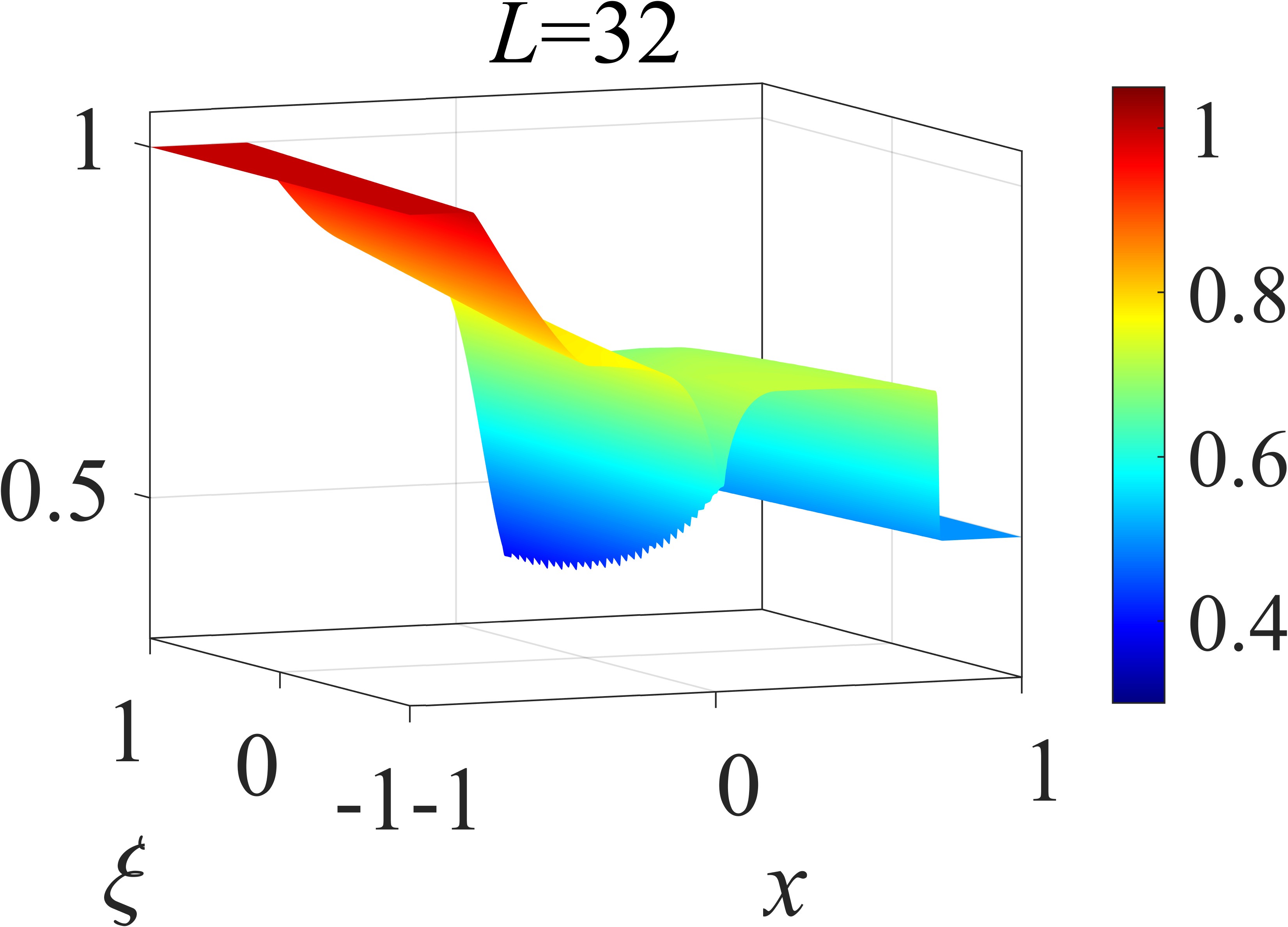}
\caption{}
\label{fig3.3c}
\end{subfigure}
\captionsetup{labelfont=bf}
\caption{Example 2: $w(x,0.8;\xi)$ obtained using (\protect\subref{fig3.3a}) gPC expansion, (\protect\subref{fig3.3b}) B-splines, and 
(\protect\subref{fig3.3c}) SP splines with different number of collocation points $L$.}
\label{fig3.3}
\end{figure}

\begin{figure}[ht!]
\centering
\begin{subfigure}[t]{0.32\textwidth}
\centering
\includegraphics[width=1.0\textwidth]{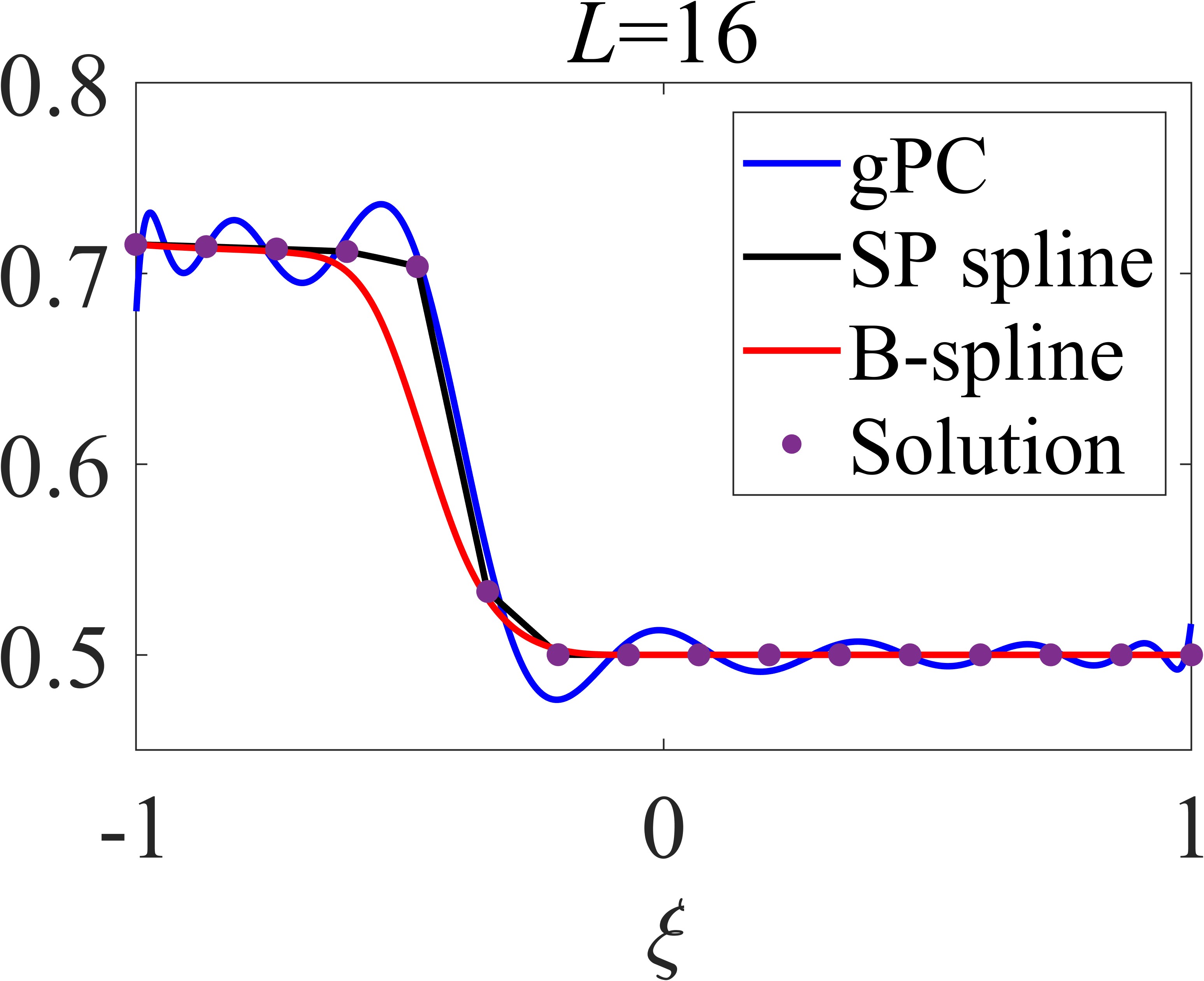}
\vspace{-0.1cm}
\end{subfigure}
\centering
\begin{subfigure}[t]{0.32\textwidth}
\centering
\includegraphics[width=1.0\textwidth]{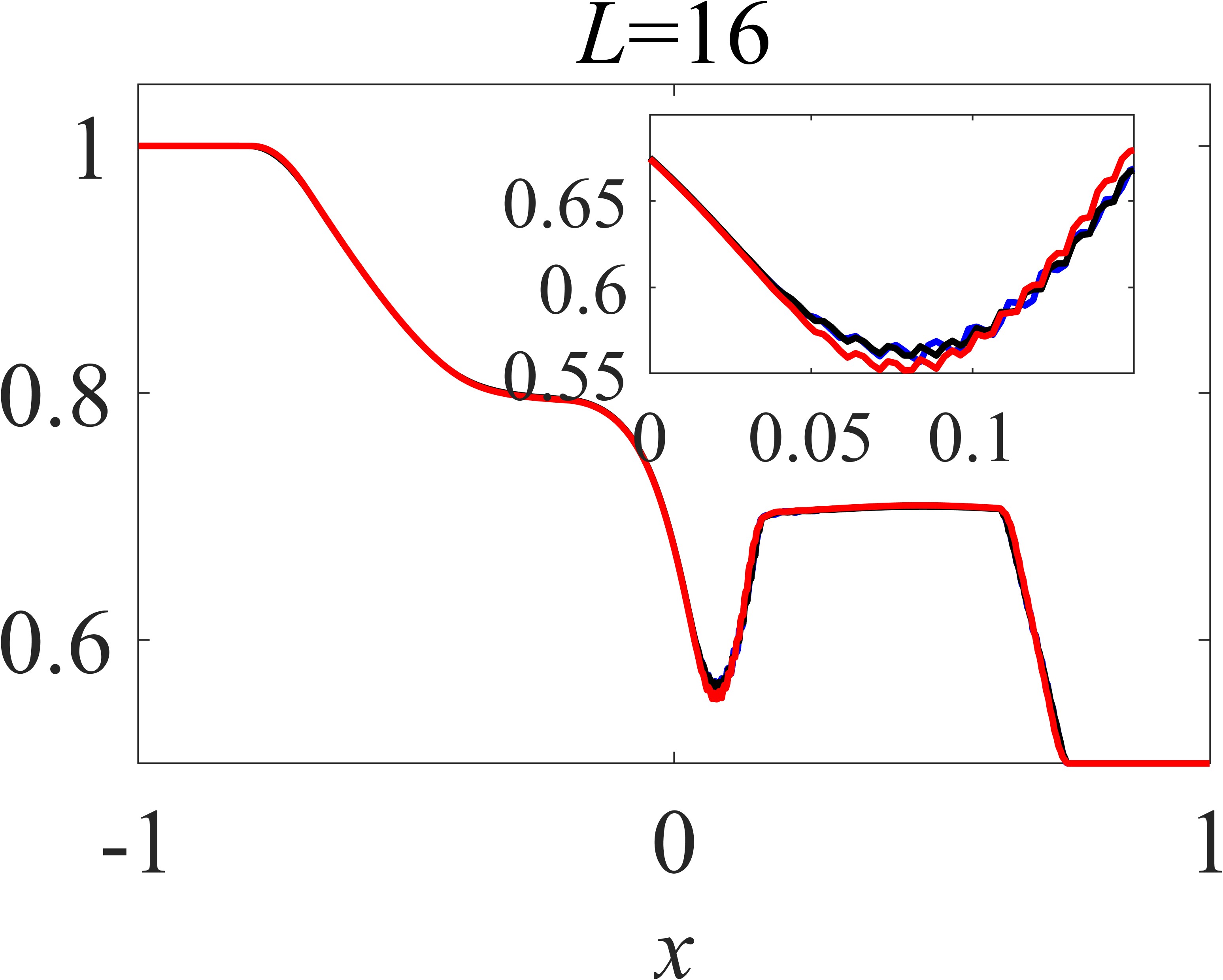}
\vspace{-0.1cm}
\end{subfigure}
\centering
\begin{subfigure}[t]{0.32\textwidth}
\centering
\includegraphics[width=1.0\textwidth]{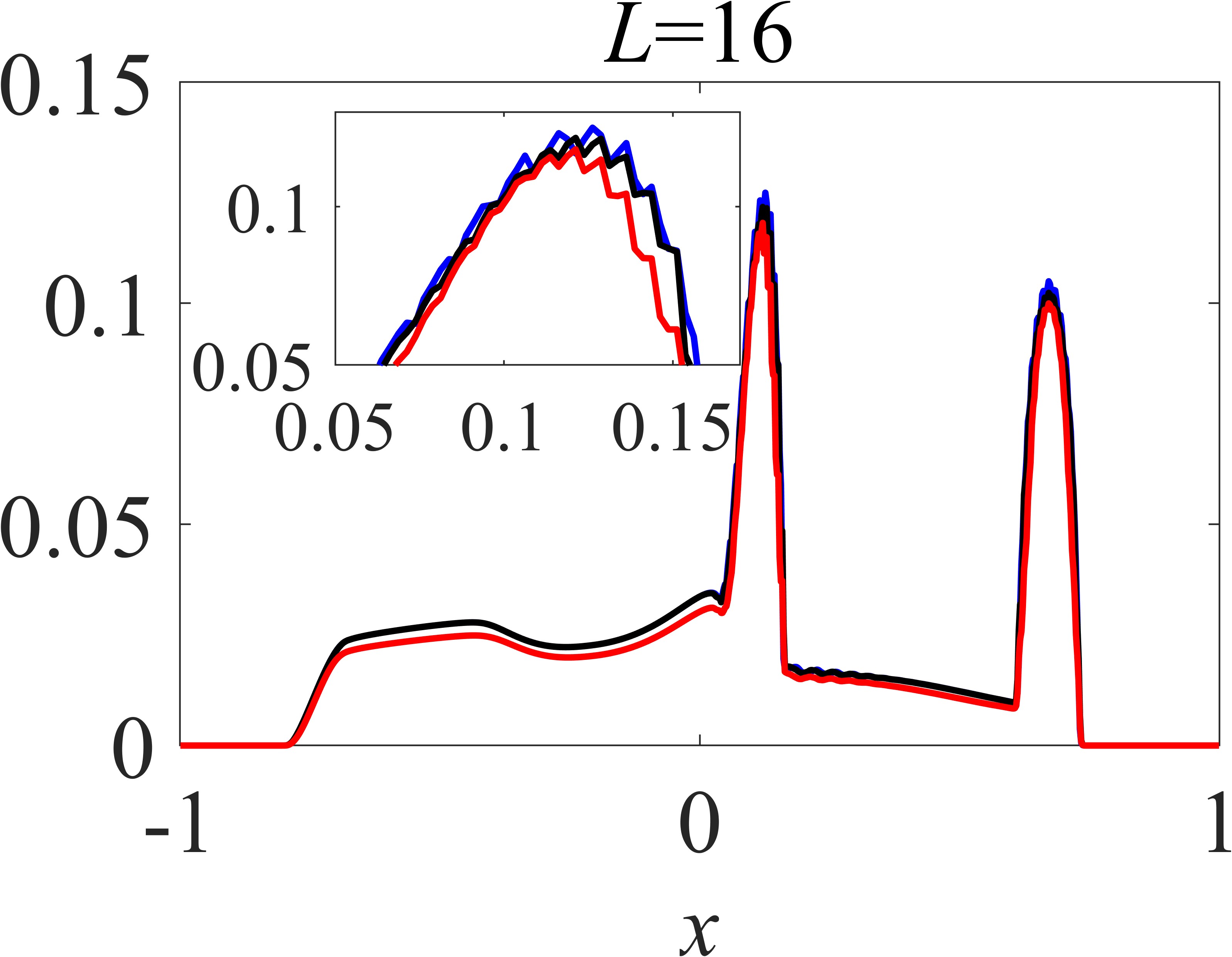}
\vspace{-0.1cm}
\end{subfigure}
\hspace{0.5cm}
\centering
\begin{subfigure}[t]{0.32\textwidth}
\centering
\includegraphics[width=1.0\textwidth]{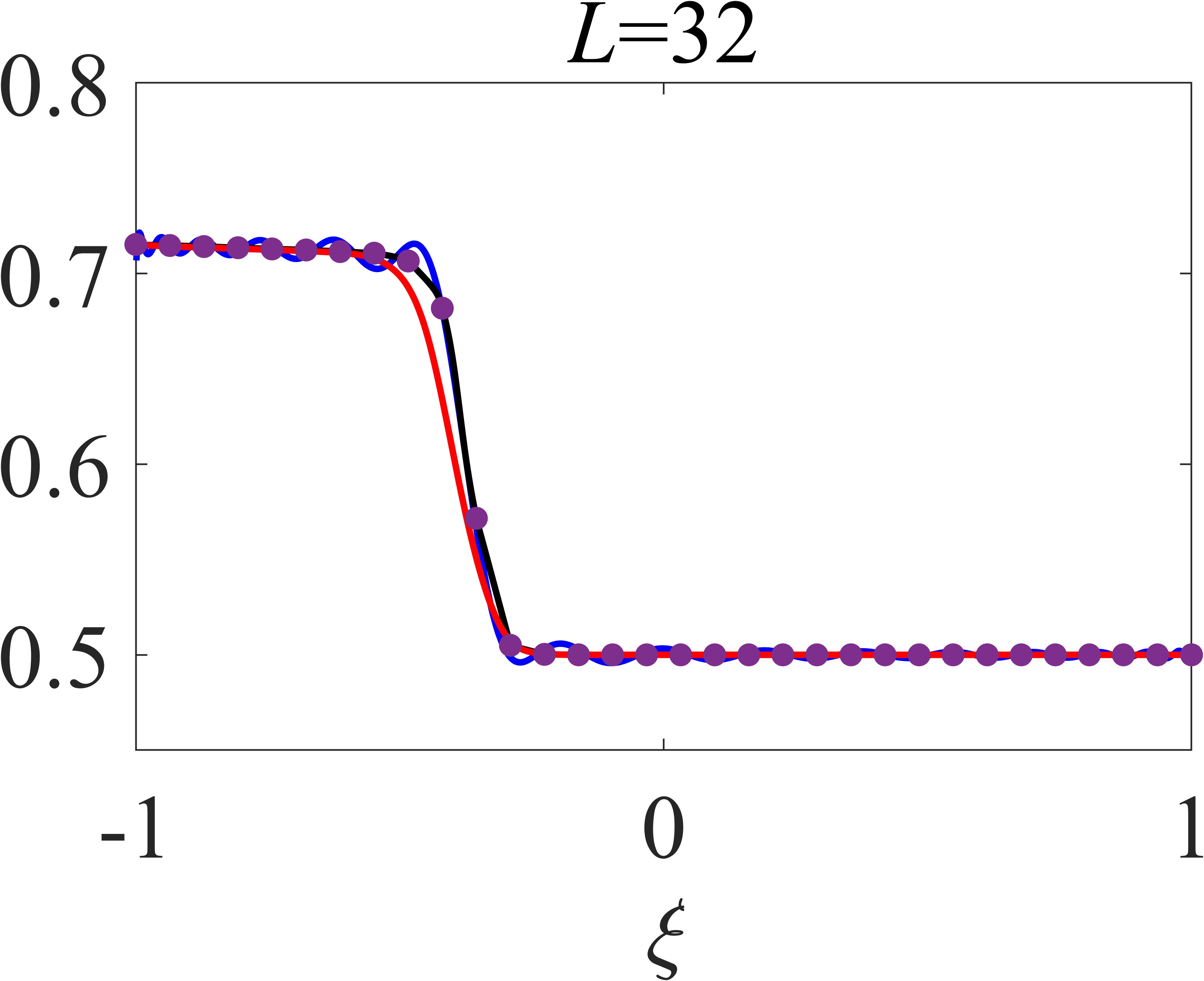}
\caption{}
\label{fig3.4a}
\vspace{-0.1cm}
\end{subfigure}
\centering
\begin{subfigure}[t]{0.32\textwidth}
\centering
\includegraphics[width=1.0\textwidth]{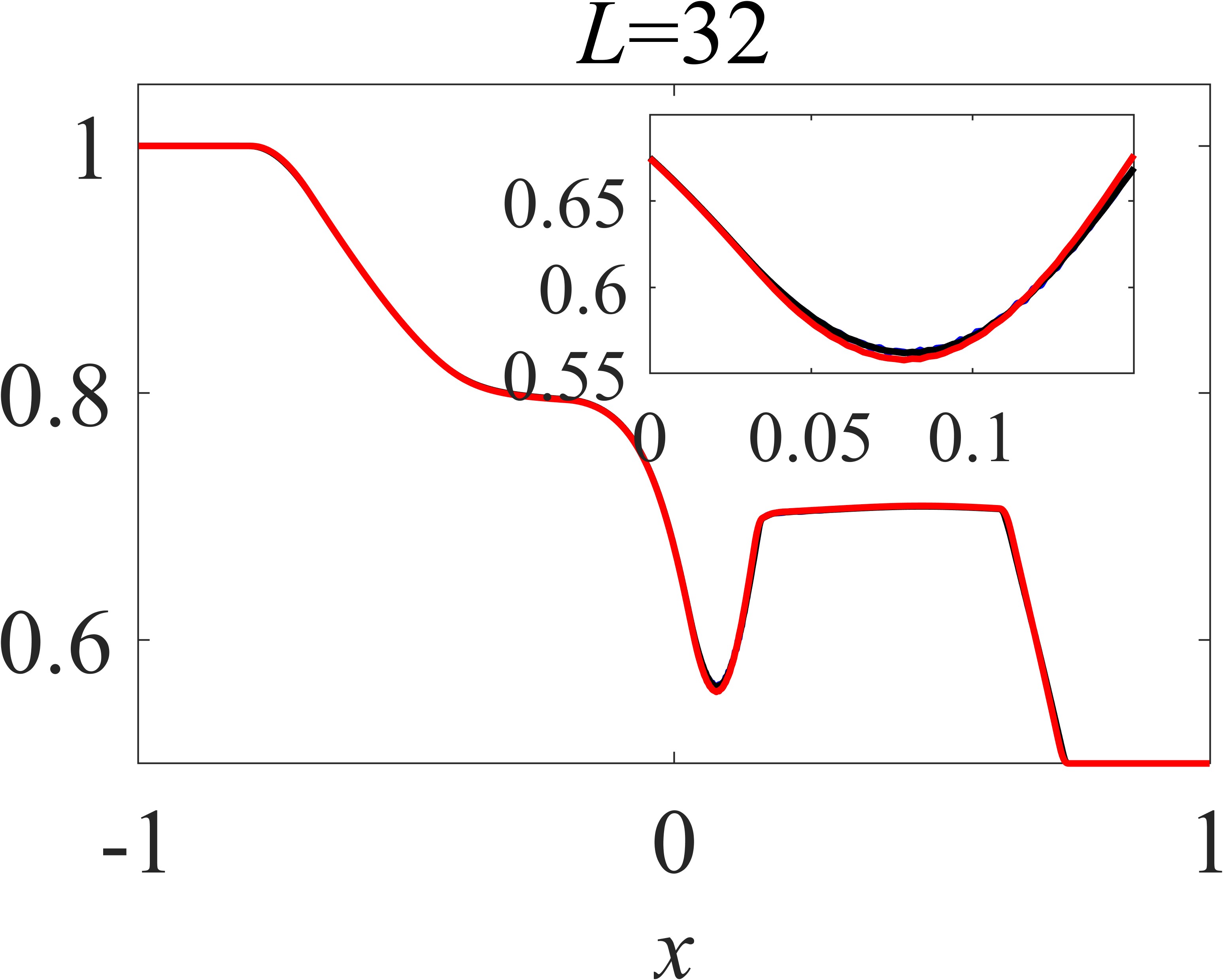}
\caption{}
\label{fig3.4b}
\vspace{-0.1cm}
\end{subfigure}
\centering
\begin{subfigure}[t]{0.32\textwidth}
\centering
\includegraphics[width=1.0\textwidth]{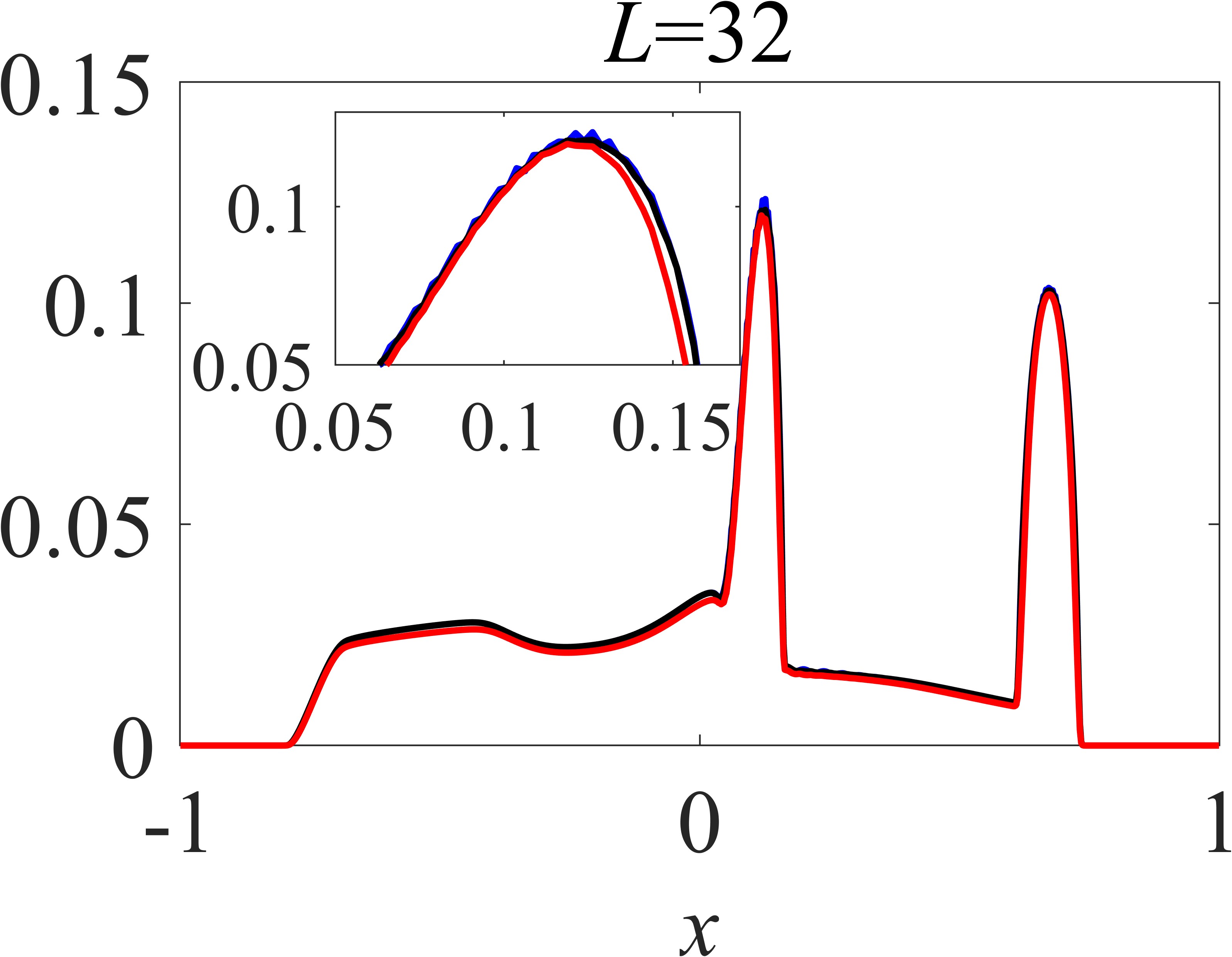}
\caption{}
\label{fig3.4c}
\vspace{-0.1cm}
\end{subfigure}
\captionsetup{labelfont=bf}
\caption{Example 2: (\protect\subref{fig3.4a}) $w(0.694,0.8;\xi)$, (\protect\subref{fig3.4b}) mean, and (\protect\subref{fig3.4c}) standard 
deviation obtained using the gPC expansion and splines. Panels in columns (b) and (c) contain a zoom at the areas of interest.}
\label{fig3.4}
\end{figure}

In \fref{fig3.3}, we plot the water surface $w(x,0.8;\xi)$. It is evident that the gPC solution exhibits oscillations near the discontinuity
location ($x\approx0.694$). These oscillations become less pronounced as the number of collocation points increases, indicating an
improvement in capturing of the stochastic behavior with increased resolution. When a spline interpolation is employed, these oscillations
are suppressed. Note that the oscillations visible near $x\approx0$ in \fref{fig3.3b} and \ref{fig3.3c} appear in the solution (not a
feature of the interpolation). Similar results are obtained for water discharge $hu$ (not shown for the sake of brevity). We additionally
plot $w(0.694,0.8;\xi)$ in \fref{fig3.4} to provide a more clear picture. The expected value and standard deviation are also shown in
\fref{fig3.4}: for all approaches, similar results are observed.

\section{Conclusions}\label{sec4}
In this paper, we have studied the stochastic collocation (SC) methods for uncertainty quantification (UQ) in hyperbolic systems of
nonlinear PDEs. We have numerically solved the underlying PDEs at a set of collocation points in random space. Then, we have used a
standard SC approach based on a gPC expansion, which relied on choosing the collocation points based on the prescribed probability
distribution and approximating the computed solution by a linear combination of orthogonal polynomials. We have illustrated that this
approach struggles to accurately capture discontinuous solutions, leading to oscillations (Gibbs-type phenomenon) that significantly deviate
from the exact solutions. We have explored alternative SC methods using uniformly distributed collocation points and employing spline
interpolations in a random space. Our study has demonstrated the effectiveness of spline-based collocation in accurately capturing and
assessing uncertainties while suppressing oscillations. We have illustrated the superiority of the spline-based collocation on two
numerical examples, including the inviscid Burgers and shallow water equations. The future work will include higher-dimensional extensions
of spline-based SC methods.

\begin{acknowledgement}
The work of A. Chertock and S. Janajra were supported in part by NSF grant DMS-2208438. The work of A. Kurganov was supported in part by
NSFC grant 12171226 and the fund of the Guangdong Provincial Key Laboratory of Computational Science and Material Design
(No. 2019B030301001). 
\end{acknowledgement}

\bibliographystyle{spmpsci}
\bibliography{ref}

\end{document}